\documentclass[11pt]{article}
\usepackage[table]{xcolor}
\usepackage{graphicx}
\usepackage{amsmath,amsfonts,amssymb}
\usepackage{graphicx}
\usepackage{multirow}
\usepackage{tikz,pgf}
\usepackage{url}
\usepackage{amsmath}
\usepackage{graphicx}
\usepackage{epsfig}
\usepackage{epstopdf}
\usepackage{color}
\usepackage{enumitem}
\def\disp{\displaystyle}

\def\tto{\;{\lower 1pt \hbox{$\rightarrow$}}\kern -10pt
\hbox{\raise 2pt \hbox{$\rightarrow$}}\;}
\def\Hat{\widehat}

\def\Bar{\overline}
\def\ra{\rangle}
\def\la{\langle}

\def\epsilon{\varepsilon}

\def\h{\hfill\Box}
\def\R{\Bbb R}

\def\N{\Bbb N}
\def\ox{\bar{x}}

\def\oy{\bar{y}}
\def\oz{\bar{z}}

\def\s{\square}

\def\ri{\mbox{\rm ri}\,}

\def\gph{\mbox{\rm gph}\,}
\def\aff{\mbox{\rm aff}\,}
\def\epi{\mbox{\rm epi}\,}

\def\dom{\mbox{\rm dom}\,}
\def\aff{\mbox{\rm aff}\,}

\def\sint{\mbox{\rm int}\,}

\def\cone{\mbox{\rm cone}}

\def\iri{\mbox{\rm iri}\,}

\newcommand{\Int}{\text{int}\,}

\def\h{\hfill\square}

\def\emp{\emptyset}

\def\oR{\Bar{\R}}

\def\gg{\gamma}

\def\bb{\beta}
\def\Th{\Theta}

\def\emp{\emptyset}

\def\oR{\Bar{\R}}

\def\gg{\gamma}

\def\Th{\Theta}

\def\qri{\mbox{\rm qri}}

\setlist[enumerate,1]{itemsep=0.0ex,parsep=0.5ex,label={\rm(\alph*)},leftmargin=*, align=left}
\newcounter{lk}

\usepackage[colorlinks=true]{hyperref}

\begin{document} \begin{center} {\sc\bf Generalized Differentiation and Duality in Infinite Dimensions\\ under Polyhedral Convexity}\\[1ex] {\sc D. V. Cuong} \footnote{Department of Mathematics,
Faculty of Natural Sciences, Duy Tan University, Da Nang, Vietnam (dvcuong@duytan.edu.vn). This research is funded by Vietnam National Foundation for Science and Technology Development (NAFOSTED)
under grant number 101.02-2020.20.}$^,$\footnote{American Degree Program, Duy Tan University, Da Nang, Vietnam.}, {\sc B. S. Mordukhovich}\footnote{Department of Mathematics, Wayne
State University, Detroit, Michigan 48202, USA (boris@math.wayne.edu). Research of this author was partly supported by the USA National Science Foundation under grants DMS-1512846 and
DMS-1808978, by the USA Air Force Office of Scientific Research grant \#15RT04, and by Australian Research Council under grant DP-190100555.}, {\sc N. M. Nam}\footnote{Fariborz Maseeh Department
of Mathematics and Statistics, Portland State University, Portland, OR 97207, USA (mnn3@pdx.edu).}, {\sc G. Sandine}\footnote{Fariborz Maseeh Department of
 Mathematics and Statistics, Portland State University, Portland, OR 97207, USA (gsandine@pdx.edu  ).}\\[2ex]
  {\bf Dedicated to Miguel Angel Goberna, in esteem}
\end{center}
\small{\bf Abstract.} This paper addresses the study and applications of polyhedral duality of locally convex topological vector (LCTV) spaces. We first revisit the classical Rockafellar's proper separation theorem for two convex sets one which is polyhedral and then present its LCTV extension replacing the relative interior by its quasi-relative interior counterpart. Then we apply this result to derive enhanced calculus rules for normals to convex sets, coderivatives of convex set-valued mappings, and subgradients of extended-real-valued functions under certain polyhedrality requirements in LCTV spaces by developing a geometric approach. We also establish in this way new results on conjugate calculus and duality in convex optimization with relaxed qualification conditions in polyhedral settings. Our developments contain significant improvements to a number of existing results obtained by Ng and Song in \cite{KFNg}. \\[1ex]

\newtheorem{Theorem}{Theorem}[section]
\newtheorem{Proposition}[Theorem]{Proposition}
\newtheorem{Remark}[Theorem]{Remark}
\newtheorem{Lemma}[Theorem]{Lemma}
\newtheorem{Corollary}[Theorem]{Corollary}
\newtheorem{Definition}[Theorem]{Definition}
\newtheorem{Example}[Theorem]{Example}
\renewcommand{\theequation}{\thesection.\arabic{equation}}
\normalsize

\section{Introduction}
\setcounter{equation}{0}

It has been well-recognized that convex separation plays a crucial role in generalized differential theory of convex analysis. One of the central results of convex and functional analysis states
that in a topological vector space, it is possible to separate a convex set with nonempty interior and another nonempty convex set with empty intersection by a closed hyperplane. Note that the
nonempty interiority condition is restrictive in both finite and infinite dimensions for applications, and thus significant effort has been made to improve this result. In his seminal monograph
``Convex Analysis" \cite{r}, Rockafellar proved that two nonempty convex sets have disjoint \emph{relative interiors} if and only if they can be separated by a hyperplane. He also proved that
this result can be further refined if one of the sets involved is a convex polyhedron. This important result and its refinement have been used broadly in convex analysis and its applications.

The important role of relative interiors in finite dimensions has been the driving force for developing generalized relative interiors notion in infinite dimensions. Among them the notion of
\emph{quasi-relative interior} introduced by Borwein and Lewis in \cite{bl} and its variants have gained certain success in developing convex generalized differentiation in infinite dimensions
with applications to convex optimization. We refer the reader to the paper by Ng and Song in \cite{KFNg} in which  the quotient topology and the notion of quasi-relative interior were used to
prove a convex separation result in LCTV spaces for a convex set and a convex polyhedron. This result and a related notion called ``strong quasi-relative
interior" were then applied to obtain the \emph{Fenchel strong duality} for convex optimization problems along with deriving generalized differential calculus rules in LCTV spaces. There is a natural
question: {\em Is the use of the strong quasi-relative interior necessary?} We will answer this question through this paper in the course of developing convex generalized differentiation for
set-valued mappings and nonsmooth functions involving convex polyhedra.

In this paper, we first revisit Rockafellar's theorem in finite dimensions as well as its extension to infinite dimensions by Ng and Song on separation of a convex set and a convex polyhedron. We
provide further details with self-contained proofs of these results for the convenience of the reader. Then we develop generalized differential calculus rules for set-valued mappings and
nonsmooth functions involving convex polyhedra with applications to the Fenchel duality theory in convex optimization. Besides obtaining new calculus rules, our developments provide significant
improvements to a number of existing results obtained by Ng and Song.

Our paper is organized as follows. In Section 2 we present basic definitions and notions of convex analysis used throughout the paper. Section 3 is devoted to revisiting the aforementioned
separation theorem by Rockafellar and its extension. We develop coderivative and subdifferential calculus rules for set-valued mappings and nonsmooth functions in Section 4. Finally, we
present applications to duality theory in Section 5.

\section{Preliminaries}
\setcounter{equation}{0}

In this section we  recall the standard notation and definitions of convex analysis in LCTV spaces used throughout the paper; see, e.g., \cite{z}.

Given a real LCTV space $X$ and its topological dual $X^*$, consider the canonical pairing $\la x^*,x\ra:=x^*(x)$ with $x\in X$ and $x^*\in X^*$. For a nonempty subset $A$ of $X$, define the
{\em conic hull} of $A$ by $\cone(A):=\big\{ta\in X\;|\;t\ge 0,\;a\in A\big\}$ and denote the {\em closure} of $A$ by $\Bar{A}$.

Recall further that two convex sets $\Omega_1,\Omega_2\subset X$ are {\em properly separated} if there exists $x^*\in X^*$ for which the following two
inequality hold:
\begin{equation}\label{ps1}
\sup\big\{\la x^*,w_1\ra\;\big|\;w_1\in\Omega_1\big\}\le\inf\big\{\la x^*,w_2\ra\;\big|\;w_2\in\Omega_2\big\},
\end{equation}
\begin{equation}\label{ps2}
\inf\big\{\la x^*,w_1\ra\;\big|\;w_1\in\Omega_1\big\}<\sup\big\{\la x^*,w_2\ra\;\big|\;w_2\in\Omega_2\big\}.
\end{equation}
Observe that condition \eqref{ps1} can be equivalently rewritten as
\begin{equation*}
\la x^*,w_1\ra\le\la x^*,w_2\ra\;\mbox{\rm whenever }w_1\in\Omega_1,\;w_2\in\Omega_2,
\end{equation*}
while \eqref{ps2} means that there exist $\bar{w}_1\in \Omega_1$ and $\bar{w}_2\in \Omega_2$ such that
\begin{equation*}
\la x^*,\bar{w}_1\ra<\la x^*,\bar{w}_2\ra.
\end{equation*}
Given two nonempty sets $\Omega_1$ and $\Omega_2$ in a topological vector space $X$. We say that $\Omega_1$ and $\Omega_2$ can be separated by a closed hyperplane that does not contain $\Omega_2$ if there exists $x^*\in X^*$ and $\alpha\in\R$ such that
\begin{equation*}
  \sup\{\la x^*,x\ra\;\big|\;x\in\Omega_1\}\leq \alpha\leq \inf\{\la x^*,x\ra\;\big|\;x\in\Omega_2\}\; \text{and}\; \alpha<\sup\{\la x^*,x\ra\;\big|\;x\in\Omega_2\}.
\end{equation*}
Given a convex subset $\Omega$ of $X$, the {\em relative interior} of $\Omega$ is defined by
\begin{equation*}\label{ri}
\mbox{\rm ri}(\Omega):=\big\{x\in\Omega\;\big|\;\exists\ \; \mbox{\rm a convex neighborhood }V\; \mbox{\rm of}\;x\;\mbox{\rm such that }\;V\cap\Bar{\aff}(\Omega)\subset\Omega\big\}.
\end{equation*}
The {\em intrinsic relative interior} of $\Omega$ is the set
\begin{equation*}
\mbox{\rm iri}(\Omega):=\big\{x\in\Omega\;\big|\;\mbox{\rm cone}(\Omega-x)\;\mbox{\rm is a subspace of }\;X\big\}.
\end{equation*}
The {\em quasi-relative interior} of $\Omega$ is the set
\begin{equation*}
\mbox{\rm qri}(\Omega):=\big\{x\in\Omega\;\big|\;\Bar{\mbox{\rm cone}(\Omega-x)}\;\mbox{\rm is a subspace of }\;X\big\}.
\end{equation*}
We say that a convex set $\Omega\subset X$ is {\sc quasi-regular} if $\qri(\Omega)={\rm iri}(\Omega)$.

Relationships between the notions of relative, intrinsic relative, and quasi-relative interiors of  convex sets in LCTV spaces give conditions that ensure quasi-regularity of a set. These results
can be found in \cite{BG,bl,hs,zduality}. \begin{Theorem}\label{ri-rel} Let $\Omega$ be a  convex subset of a topological vector space $X$. Then we have the inclusions
\begin{equation}\label{ri-rel1} \ri(\Omega)\subset\iri(\Omega)\subset\qri(\Omega). \end{equation} If furthermore $X$ is locally convex and  $\ri(\Omega)\ne\emp$, then the inclusions in
\eqref{ri-rel1} become the equalities \begin{equation*}\label{eq-ri-quari} \ri(\Omega)=\iri(\Omega)=\qri(\Omega). \end{equation*} \end{Theorem}

We continue with a well-known version of {\em proper separation} of a singleton from a convex set that gives us yet another characterization of quasi-relative interior; see: \cite[Theorem 2.3]{Flores1}.

\begin{Proposition}\label{qri-sep1} Let $\Omega$ be a convex set in an LCTV space $X$, and let $\ox\in\Omega$. Then the sets $\{\ox\}$ and $\Omega$ are properly separated if and only if
$\ox\notin\qri(\Omega)$.
\end{Proposition}

 A subset $P\subset X$ is said to be a  polyhedral set if there exist $x_i^*\in X^*,\alpha_i\in\mathbb{R}$ for $i=1,\ldots,k,$  such that
\begin{equation*}\label{Poly}
  P=\{x\in X\ | \  \la x_i^*,x\ra\leq \alpha_i,i=1,\ldots, k\}.
\end{equation*}

Let $\Omega$ be a nonempty convex subset of $X$. The {\sc normal cone} to $\Omega$ at $\ox$ is defined by
\begin{equation}\label{nc}
N(\ox;\Omega):=\big\{x^*\in X^*\;\big|\;\la x^*,x-\ox\ra\le 0\;\text{ for all }\;x\in\Omega\big\}\; \text{whenever}\; \ox\in\Omega
\end{equation}
and $N(\ox;\Omega):=\emp$ if $\ox\notin\Omega$.

A set-valued mapping $F\colon X\tto Y$ between topological vector spaces $X$ and $Y$ is said to be {\em convex} if its {\em graph} given by
\begin{equation*}
    \gph(F):=\{(x, y)\; |\; y\in F(x)\}
\end{equation*}
is a convex set in $X\times Y$. The {\em domain} of $F$ is the set of all $x\in X$ for which $F(x)$ is nonempty.

Given $(\ox, \oy)\in \gph(F)$, the {\em coderivative} of $F$ at $(\ox, \oy)$ is the set-valued mapping $D^*F(\ox, \oy)\colon Y^*\tto X^*$ defined by
\begin{equation}\label{cod}
    D^*F(\ox, \oy)(y^*):=\{x^*\in X^*\; |\; (x^*, -y^*)\in N((\ox, \oy); \gph(F)\}, \; y^*\in Y^*.
    \end{equation}
    Note that if $A\colon X\to Y$ is a continuous linear mapping, then
    \begin{equation*}
        D^*A(\ox, \oy)(y^*)=\{A^*(y^*)\}, \; y^*\in Y^*,
    \end{equation*}
    where $\oy:=A(\ox)$ and $A^*$ denotes the adjoint mapping of $A$.

    Given a convex function $f\colon X\to \oR:=(-\infty, \infty]$, the {\em subdifferential} of $f$ at $\ox\in\dom(f):=\{x\in X\; |\; f(x)<\infty\}$ is defined by
    \begin{equation*}
        \partial f(\ox):=\{x^*\in X^*\; |\; \la x^*, x-\ox\ra\leq f(x)-f(\ox)\; \mbox{\rm for all }x\in X\}.
    \end{equation*}
    Defining the {\em epigraphical mapping} $E_f\colon \R\tto \R$ by $E_f(x):=[f(x), \infty)$, we have the following representation
    \begin{equation*}
        \partial f(\ox)=D^*E_f(\ox, \oy)(1), \; \mbox{\rm where }\oy:=f(\ox).
    \end{equation*}

\section{Rockafellar's Polyhedral Separation Theorem Revisited} \setcounter{equation}{0}

The following theorem was known by \cite[Theorem 20.2]{r} or \cite[Proposition 1.5.7]{Bertsekas2003}. We present here a rather different and more detailed proof.

\begin{Theorem}{\bf(Rockafellar's theorem for polyhedral sets).}\label{TheoRoc-Po} Let $P$ and $\Omega$ be two nonempty convex sets of $\mathbb{R}^n$ such that $P$ is polyhedral. The condition
that there exists a (closed) hyperplane  separating $P$  and $\Omega$ properly and not containing $\Omega$ is equivalent to \begin{equation*} P \cap \ri(\Omega)=\emp. \end{equation*} \end{Theorem} {\bf
Proof.} Assume there exists a hyperplane $H$ defined by \begin{equation*}
  H:=\{x\in \mathbb R^n\; |\; \la v,x\ra =\alpha\}\text{ for some } v\in \mathbb{R}^n,\,v\neq 0, \text{ and } \alpha\in\mathbb{R}
\end{equation*}
that separates $P$  and $\Omega$ properly and does not contain $\Omega.$ Since $H$ properly separates $P$  and $\Omega$, without loss of generality we have
\begin{equation}\label{eq-sepa-O-P}
  \la v,x\ra \leq \alpha \leq \la v,p\ra \text{ for all}\ x\in\Omega,\,p\in P,
\end{equation}
and since $\Omega\not\subset H$ there exists $\ox\in \Omega$ such that
\begin{equation*}
\la v,\ox\ra < \alpha \leq \la v,p\ra \text{ for all } p\in P.
\end{equation*}
Next, observe that
\begin{equation*}
 \la v,x\ra \leq \alpha=\la v, u\ra \text{ for all }x\in \Omega, u\in H
\end{equation*}
and that the element $\ox$ satisfies
\begin{equation*}
 \la v,\ox\ra < \alpha=\la v, u\ra \text{ for all }u\in H,
\end{equation*} therefore the hyperplane $H$ properly separates the convex sets $H$ and $\Omega$. To arrive at a contradiction, suppose that \begin{equation*}
  P \cap \ri(\Omega) \ne\emp
\end{equation*}
and choose an element $x_0\in  P \cap \ri(\Omega)$. Then by \eqref{eq-sepa-O-P} we have
\begin{equation*}
 \la v,x_0\ra \leq \alpha \leq \la v, x_0\ra,
\end{equation*}
or $\la v,x_0\ra =\alpha$. This implies $x_0\in H$ but $\ri(H)=H$ since $H$ is an affine set, therefore
\begin{equation*}
x_0 \in H \cap \ri(\Omega) =\ri(H) \cap  \ri(\Omega)
\end{equation*}
which contradicts that $H$  and $\Omega$ are properly separated since their proper separation is equivalent to $\ri(H)\cap \ri(\Omega)=\emp$. Therefore, $P\cap \ri(\Omega)=\emp$.

Conversely, suppose that $P\cap \ri(\Omega)=\emptyset$. We will show that there exists a hyperplane  that properly separates $\Omega$ and $P$, while it does not contain $\Omega$. Define the set
\begin{equation*} D:=P \cap \aff(\Omega). \end{equation*} Since $P$ and $\aff(\Omega)$ are convex polyhedra, if $D=\emptyset,$ then the strict separation theorem gives that $P$ and $\aff(\Omega)$
can be strictly separated by a hyperplane $H$. This hyperplane strictly separates $P$ and $\Omega$ and does not contain $\Omega$ since $\Omega\subset\aff(\Omega)\subset\Int(H^+)$ while
$P\subset\Int(H^-)$.

We now consider the case where $D\ne\emptyset$. It follows from the assumption $P \cap \ri(\Omega)=\emptyset$ that
\begin{equation*}
 \ri(D)\cap \ri(\Omega)\subset D\cap \ri(\Omega)\subset P\cap \ri(\Omega)=\emp.
\end{equation*} By the proper separation theorem, the sets $D$ and $\Omega$ can be properly separated by a hyperplane $H_0$ defined by
  \begin{equation*}
 H_0:=\{x\in \mathbb R^n\ |\ \la v_0,x\ra =\alpha_0\}\ \text{for some}\ v_0\in \mathbb{R}^n, v_0\neq 0,  \text{and}\ \alpha_0\in\mathbb{R}.
\end{equation*}
If $\Omega\subset H_0,$ then $\aff(\Omega)\subset H_0$ and hence $D\subset H_0$ therefore $\la v_0,x\ra=\la v_0,y\ra=\alpha_0$ for all $x\in\Omega$, $y\in D$. This contradicts the proper separation of
$D$ and $\Omega$ by $H_0$, so $H_0$ does not contain $\Omega$.

Suppose that \begin{equation*}
  \Omega\subset H_0^-:=\{x\in \mathbb R^n\ |\ \la v_0,x\ra \leq \alpha_0\}\textrm{ so }D\subset H_0^+:=\{x\in \mathbb R^n\ |\ \la v_0,x\ra \geq \alpha_0\},
\end{equation*}
and define
\begin{equation*}
 \Theta:=H_0^-\cap \aff(\Omega).
\end{equation*}
Since $H_0$ does not contain $\Omega$ and $\Omega\subset\Theta$, $H_0$ also does not contain $\Theta,$ and hence the sets $H_0$ and $\Theta $ are properly separated, thus
\begin{equation*}
 H_0\cap \ri(\Theta)=\ri(H_0)\cap \ri(\Theta)=\emptyset.
\end{equation*}
Also, if $x\in D\cap\ri(\Theta)$ then $x\in\Theta$ gives $\la v_0,x\ra\leq\alpha_0$ and $x\in D$ gives $\alpha_0\leq\la v_0,x\ra$, or $\la v_0,x\ra=\alpha_0$ so $x\in H_0$. This contradicts
$H_0\cap\ri(\Theta)=\emp$, thus $D\cap\ri(\Theta)=\emp$. Additionally, if $x\in P\cap\ri(\Theta)$ then $x\in P\cap\aff(\Omega)=D$ since $\ri(\Theta)\subset\aff(\Omega)$, which means $x\in D\cap\ri(\Theta)$ contradicting $D\cap\ri(\Theta)=\emp$. Thus, we also have
\begin{equation}\label{eq-Pcap-The}
P\cap\ri(\Theta)=\emp.
\end{equation}

Now if $P\cap\Theta=\emp$, since $P$ and $\Theta$ are polyhedral, by the strict separation theorem there exists a hyperplane $H$ strictly separating $P$ and $\Theta$ such that \begin{equation*}
\Omega\subset\Theta\subset\Int(H^-)\textrm{ and }P\subset\Int(H^+). \end{equation*} This hyperplane $H$ separates $P$ and $\Omega$ and does not contain $\Omega$.

We now consider the final case where $P\cap \Theta\ne\emptyset$. We can suppose that $0\in P\cap \Theta$, and since $\Theta\subset\aff(\Omega)$ we have
\begin{equation*}
0\in P\cap\Theta\subset P\cap\aff(\Omega)=D\textrm{ which gives } 0\in\Theta\cap D\subset H_0^-\cap H_0^+=H_0.
\end{equation*}
Since $0\in H_0$, we have $\alpha_0=0$ or $H_0=\{x\in\R^n:\la v_0,x\ra=0\}$.

Observe that $0\notin\Int(P)$. Indeed, to arrive at a contradiction, suppose that $0\in\Int(P)$. For each $\ox\in\ri(\Theta)$ since $0\in\Theta$ we have $[\ox,0)\subset \ri(\Theta)$. But
$0\in\Int(P)$ would mean $\Int(P)\cap[\ox,0)\ne\emptyset$ so $\Int(P)\cap\ri(\Theta)\not=\emp$ contradicting that $P\cap\ri(\Theta)=\emptyset$, and hence $0\notin\Int(P)$.

Since $0\in P\setminus{\rm int}(P)$ and $P$ is a convex polyhedron, $P$ can be represented as the intersection of half-spaces
\begin{equation*}   P=\{x\in \mathbb{R}^n\ |\ \la u_i,x\ra \leq 0,\ i=1,\ldots, m\}\cap \{x\in \mathbb{R}^n\ |\ \la u_j,x\ra \leq \beta_j,\ j=m+1,\ldots, \bar{m}\}
\end{equation*}
where $1\leq m\leq\bar{m}$, $v_k$ are nonzero vectors in $\mathbb{R}^n$ for $k=1,\ldots,\bar{m},$ and $\beta_j>0$ for $j=m+1,\ldots,\bar m$.

Let $M$ be the relative boundary of $\Theta,$ i.e., \begin{equation*} M=H_0\cap \aff(\Omega) \end{equation*} as in Figure 1, and note that $M$ and $\aff(\Omega)$ are both
linear subspaces due to $0\in H_0$ and $0\in\Theta\subset\aff(\Omega)$, so they are both affine sets containing the origin. Additionally, the half-subspace $\Theta$ is a cone.

\begin{figure}[h!bt]\label{rockthmP}
\centering
\includegraphics[width=4in]{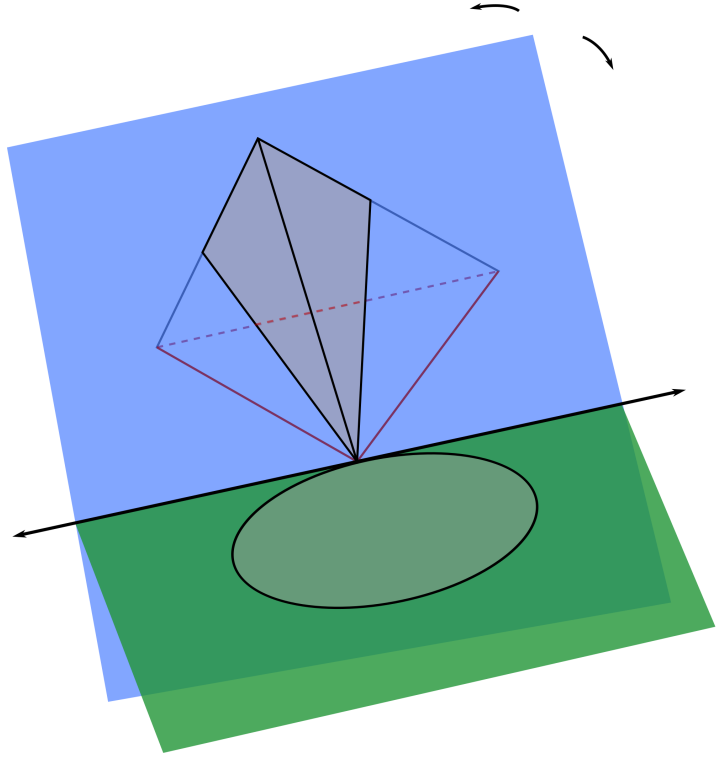}
\caption{$\Omega$, $P$, $D$, $H_0$, $M$, and $\Theta$}
\end{figure}

Now, consider the polyhedral convex cone \begin{equation*} K:=\{x\in \mathbb{R}^n\ |\ \la u_i,x\ra \leq 0,\ i=1,\ldots, m\} +M = \cone(P)+M. \end{equation*} We claim that $K\cap
\ri(\Theta)=\emptyset.$ To arrive at a contradiction, suppose that there exists $\ox\in K\cap\ri(\Theta).$ Since $\ox\in K$ and $\ox\notin M$, \begin{equation*}
  \ox=\bar{\gamma}\bar{w}+\bar{u},\ \text{for some }\ \bar{w}\in P,\bar{u}\in M,\bar{\gamma}>0.
\end{equation*}
 It follows that \begin{equation*} \bar{w}=\frac{1}{\bar{\gamma}}\ox-\frac{1}{\bar{\gamma}}\bar{u}\in P. \end{equation*} Since $\ox\in \ri(\Theta)$ and since $-\bar{u}\in
M\subset\Theta$ because $M$ is a subspace, one has \begin{equation}\label{eq-ri-lin}
  \lambda \bar{x}-(1-\lambda)\bar{u}\in \ri(\Theta)\textrm{ for all }\lambda \in(0,1).
\end{equation} Since $\Theta$ is a cone, $\ri(\Theta)$ is a cone as well. So \eqref{eq-ri-lin} together with the fact that $\ri(\Theta)$ is a cone gives us the inclusion
 \begin{equation*}
     \frac{\lambda}{1-\lambda} \bar{x}-\bar{u}\in \ri(\Theta) \textrm{ for all } \lambda \in(0,1)
 \end{equation*}
 which is equivalent to
 \begin{equation*}
    \alpha \bar{x}-\bar{u}\in \ri(\Theta) \textrm{ for all } \alpha >0.
\end{equation*}
In particular, we have $\bar x-\bar u\in\ri(\Theta)$, thus multiplying by $\frac{1}{\bar\gamma}>0$ gives
 \begin{equation*}
 \bar{w}=\frac{1}{\bar{\gamma}}\ox-\frac{1}{\bar{\gamma}}\bar{u}\in \ri(\Theta)
\end{equation*} which means $\bar{w}\in P\cap\ri(\Theta).$ This contradicts \eqref{eq-Pcap-The}, hence $K\cap\ri(\Theta)=\emptyset.$

Since $K$ is a convex polyhedra, $K$ can be represented as
\begin{equation*}
  K=\bigcap_{i=1}^kH_i^-,
\end{equation*}
where for $i=1,\ldots,k$, each $H_i^-$ is a closed half-space with respect to a hyperplane $H_i$ which passes through the origin since $K$ is a cone. That is, for each $i$, there is a nonzero vector
$v_i\in\R^n$ such that $H_i^-=\{x\in\R^n:\la v_i,x\ra\leq 0\}$. Additionally, since $M\subset K$ we have $M\subset H_i^-$ for all $i$.

We are going to show that $H_i^-\cap \ri(\Theta)=\emp$ for at least one $i$. Note first that if $H_i^-\cap\ri(\Theta)\neq\emp$, then $\Theta\subset H_i^-$. Indeed, suppose that $\bar x\in
H_i^-\cap\ri(\Theta)$ so
 $\la v_0,\bar x\ra<0$, choose $y\in\Theta$, and let
\begin{equation*} \lambda=\frac{\la v_0,y\ra}{\la v_0,\bar x\ra}\textrm{ and }z=y-\lambda\bar x. \end{equation*} Now $\lambda\geq 0$, $z\in H_0$ since $\la v_0,z\ra=0$, and $z\in\aff(\Omega)$
since $\bar x,y\in\aff(\Omega)$ and $\aff(\Omega)$ is a subspace. This gives $z\in H_0\cap\aff(\Omega)=M$, therefore $z\in H_i^-$ as well. Since $y=\lambda\bar x+z$, $\la v_i,z\ra\leq 0$, and
$\la v_i,\bar x\ra\leq 0$, it follows that \begin{equation*} \la v_i,y\ra=\la v_i,\lambda\bar x+z\ra=\lambda\la v_i,\bar x\ra+\la v_i,z\ra\leq 0 \end{equation*} therefore $y\in H_i^-$. Since $y$
represents an arbitrary element of $\Theta$, we conclude $\Theta\subset H_i^-$ whenever $H_i^-\cap\ri(\Theta)\neq\emp$. By this result, it must be the case that $H_{i_0}\cap\ri(\Theta)=\emp$ for
some $i_0\in\{1,\ldots,m\}$. If not, \begin{equation*} \ri(\Theta)\subset\Theta\subset\cap_{i=1}^k H_i^-=K \end{equation*} which contradicts $\ri(\Theta)\cap K=\emp$. If we let $H=H_{i_0}$ for
some index $i_0$ where $H_{i_0}^-$ has an empty intersection with $\ri(\Theta)$, we have $\Omega\subset\Theta\subset H^+$, $P\subset K\subset H^-$, and $\Omega\not\subset H$ since
$\ri(\Theta)\subset\Theta\not\subset H$. $\h$

Let  $P$ be a general polyhedral set. It follows from \cite[Proposition 2.197]{Bonnans2000} that $\ri(P)\ne\emptyset.$ Therefore, \cite[Theorem 2.12]{BG} shows that $\ri(P)=\iri(P)
=\qri(P).$ This means that a convex polyhedron is quasi-regular.

Next we discuss an extension of Theorem \ref{TheoRoc-Po} to locally convex topological vector spaces obtained by Ng and Song; see \cite[Theorem 3.1]{KFNg}. Their proof uses the quotient topology
to reduce the general result to the finite-dimensional case. Here are some elaborations in this direction.

Given a linear subspace $N$ of a topological vector space $X$, recall that the \emph{quotient space} $X/N$ is defined by
\begin{equation*}
X/N:=\big\{x+N\; \big |\; x\in X\big\}.
\end{equation*}
The addition and the scalar multiplication on $X/N$ are defined by
\begin{equation*}
(x+N)+(y+N):=(x+y)+N\; \mbox{\rm and }\alpha (x+N):=\alpha x+N
\end{equation*}
for $x, y\in X$ and scalar $\alpha$. Since $N$ is a linear subspace, both operations above are well-defined. We can show easily that $X/N$ with these operations is a vector space.

The quotient map $\pi\colon X\to X/N$ is defined by
\begin{equation}\label{pi} \pi(x):=x+N\; \mbox{\rm for }x\in X. \end{equation} Now, we denote by $\tau_N$ the strongest topology on $X/N$ such that the
quotient map is continuous, i.e., \begin{equation*} \tau_N=\{E\subset X/N\; \big|\; \pi^{-1}(E)\; \mbox{\rm is open in }X\}. \end{equation*}

The following theorem can be deduced from the definitions; see \cite{Rudin1991}.

\begin{Theorem} Let $X$ be a topological vector space, and let $N$ be a closed subspace of $X$. Then we have the following assertions:
\begin{enumerate}
\item The quotient map $\pi\colon X\to X/N$ is linear, continuous, and open.
\item The vector space $X/N$ with the quotient topology is a topological vector space.
\item If $\mathcal{B}$ is a basis of neighborhoods of the origin in $X$, then $\pi(\mathcal{B})$ is a basis of neighborhoods of the origin in $X/N$.
\item If $X$ is an LCTV space, then so is $X/N$.
\end{enumerate}
\end{Theorem}

Next we present some useful propositions of their own interest.

\begin{Proposition}\label{WS1} Let $X$ be a topological vector space, and let $N$ be a closed linear subspace of $X$. Then the adjoint $\pi^*\colon (X/N)^*\to X^*$ of $\pi\colon X\to X/N$ is a
bijection from $(X/N)^*$ to $N^\perp:=\{f\in X^*\; |\; f(x)=0\; \mbox{\rm for all }x\in N\}$.
\end{Proposition}
{\bf Proof.}  Fix any $z^*\in (X/N)^*$. We will show that $\pi^*(z^*)\in N^\perp$. For any $x\in N$ we have $\pi(x)=0$, the zero element in $X/N$, and so
\begin{equation*}
\la \pi^*(z^*), x\ra =\la z^*, \pi(x)\ra =0.
\end{equation*}
It follows by the definition that $\pi^*(z^*)\in N^\perp$.

To prove the reverse inclusion, we fix any $f\in N^\perp$ and show that $f=\pi^*(z^*)$ for some $z^*\in (X/N)^*$. Now, fix any $z:=x+N\in X/N$ for $x\in X$. Define \begin{equation*} z^*(z):=f(x).
\end{equation*} Since $f\in N^\perp$, we can check that $z^*$ is well-defined. Obviously, $z^*$ is linear, and its continuity follows from the continuity of $f$ and the fact that $\pi$ is an open
map. Thus $z^*\in
(X/N)^*$ and \begin{equation*} \pi^*(z^*)=z^*\circ \pi =f, \end{equation*} which shows that $\pi^*$ is surjective.

Finally, suppose that $z^*\in(X/N)^*$ and $\pi^*(z^*)=0$. In that case, for any $x\in X$ we have $0=\la\pi^*(z^*),x\ra=\la z^*,\pi(x)\ra$ which means $z^*=0$, thus $\pi^*$ is injective as well.
This completes the proof. $\h$

\begin{Proposition}\label{WS2} Let $X$ be an LCTV space, and let $f_i\colon X\to \Bbb F$ be  continuous linear functionals for $i=1, \ldots, m$ with $N:=\cap_{i=1}^m \mbox{\rm ker }(f_i)$. Then $N$ is a
closed linear subspace of $X$ and $\mbox{\rm dim}(X/N)<\infty$.
\end{Proposition}
{\bf Proof.} Obviously, $N$ is a closed linear subspace of $X$. Define the map $F\colon X\to \Bbb F^m$ by $F(x):=(f_1(x), f_2(x), \ldots, f_m(x))$. Then $F$ is a linear mapping with $\mbox{\rm ker}(F)=N$
and thus
\begin{equation*}
\mbox{\rm dim}(X/N)=\mbox{\rm dim}(F(X))\leq m,
\end{equation*}
which completes the proof. $\h$

We are now ready to present  the afore-mentioned  extension of Theorem \ref{TheoRoc-Po} to locally convex topological vector spaces obtained by Ng and Song (\cite[Theorem 3.1]{KFNg}).

\begin{Theorem}\label{Sepa-poly-convex} Let $P$ and $\Omega$ be nonempty convex sets in an LCTV space $X.$ Suppose that $P$ is a convex polyhedron and $\qri(\Omega)\ne\emptyset.$ Then
$P \cap \qri(\Omega)=  \emptyset$ if and only if there exists a closed hyperplane  $H$ such that $H$ does not contain $\Omega$ and separates $P$ and $\Omega$ properly. \end{Theorem} {\bf Proof.}
Suppose that the sets $P$  and $\Omega$ can be properly separated by a closed hyperplane $H$ and $\Omega\not\subset H.$ Then there exist $f\in X^*$ and $\alpha\in \R$ such that \begin{equation}\label{eqTheo-p-in1}
  \sup_{x\in P}f(x)\leq\alpha \leq \inf_{x\in \Omega} f(x) \; \mbox{\rm and }\alpha <\sup_{x\in \Omega} f(x).
\end{equation}
On a contrary suppose that  $P\cap \qri(\Omega)\ne \emptyset.$ Then there exists $\ox\in P\cap \qri(\Omega)$. It follows from \eqref{eqTheo-p-in1} that $f(x) \geq 0$ for all  $x\in \overline{\cone}
(\Omega-\ox).$ Since $\overline{\cone}(\Omega-\ox)$ is a linear subspace, $f(x)=0$ for all $x\in\Omega-\ox,$ which is a contradiction to \eqref{eqTheo-p-in1}.

For the converse statement, we assume that $P\cap \qri(\Omega)=\emptyset$ with $\qri(\Omega)\ne\emptyset$ and show the fulfillment of the separation property formulated
in the theorem. Since
$P$ is a convex polyhedron, we have the following representation
\begin{equation*}
  P=\{x\in X\; |\; f_i(x) \leq b_i,\ i=1,\ldots,m\},
\end{equation*}
where $f_i\in X^*, b_i\in\R$ for $i=1, \ldots, m$. Set $N:=\cap_{i=1}^m{\rm ker}(f_i).$ By Proposition \ref{WS2}, the quotient space $X/N$ is finite-dimensional. Let $\pi \colon X \to  X/N$ be the
quotient map. For each $i\in\{1, \ldots, m\}$ consider the functions $\Hat{f}_i\colon X/N\to \R$ given by \begin{equation*}
  \Hat{f}_i([x]):=f_i(x),\; \text{where}\; [x]:=x+N\in X/N.
\end{equation*} It is not hard to show that $\Hat{f}_i$ is well-defined and  $\Hat{f}_i\in (X/N)^*$ for all $i=1, \ldots, m$.  Remembering the construction of the quotient map $\pi\colon X\to X/N$ from \eqref{pi} gives us easily that
\begin{equation*}
  \pi(P)=\{[x]\in X/N\; |\; \Hat{f}_i([x]) \leq b_i,\ i=1,\ldots,m\},
\end{equation*}
and hence $\pi(P)$ is a convex polyhedron in $X/N.$ Since  $X/N$ is finite-dimensional, we have $\pi\big(\qri(\Omega)\big)=\ri\big(\pi(\Omega)\big).$ Assuming now that  $\pi(P)\cap \ri (\pi(\Omega))
\ne \emptyset$, i.e., there exists $\ox\in\qri(\Omega)$ with $[\ox]\in \pi(P)$, implies that
\begin{equation*}
  f_i(\ox)= \Hat{f}_i([\ox])\ra \leq b_i\; \text{for all}\; i=1, \ldots, m,
\end{equation*}
and $\ox\in P\cap\qri(\Omega)$. This contradicts our initial assumption. Thus $\pi(P)\cap \ri\big(\pi(\Omega)\big)=\emptyset.$  By Theorem \ref{TheoRoc-Po}, there exists $x^*\in (X/N)^*$ and $\alpha\in\R$ such
that
\begin{equation*}
  \sup_{[x]\in \pi(P)}\la {x}^*,[x]\ra \leq \alpha\leq \inf_{[x]\in \pi(\Omega)}\la {x}^*,[x]\ra\;\ \text{and}\; \alpha< \sup_{[x]\in \pi(\Omega)}\la {x}^*,[x]\ra.
\end{equation*} Proposition \ref{WS1} ensures that the adjoint operation $\pi^*\colon (X/N)^*\to X^*$ is into the space $N^\perp$. Letting $f:=\pi^*(x^*)$, we deduce from this proposition that $f\in N^\perp\subset X^*$ with $f(x)=\la
x^*, [x]\ra$ for $x\in X$, and so \begin{equation*} \sup_{x\in P} f(x) \leq \alpha\leq \inf_{x\in \Omega}f(x)\; \text{and}\;  \alpha <\sup_{x\in \Omega} f(x), \end{equation*} which completes the proof of the theorem.
$\h$

\begin{Remark} {\rm In our proof of Theorem \ref{Sepa-poly-convex}, we corrected a glitch in the proof of \cite[Theorem 3.1]{KFNg} by Ng and Song claiming that if $\Omega$ and $P$ can be properly separated by a closed hyperplane $H$ and $\Omega\not\subset H$, then there exists $f\in X^*$  such that \begin{equation*}
  \sup_{x\in P}f(x) \leq \inf_{x\in \Omega} f(x)<\sup_{x\in \Omega} f(x).
\end{equation*}
Note that this is not the case in general unless $P\cap \Omega\neq \emptyset$.}
\end{Remark}

\section{Polyhedral Normal Cone Intersection Rule} \setcounter{equation}{0}

The following normal cone intersection rule is the driving force of our {\em geometric approach} to generalized differential calculus developed in what follows; cf.\ \cite{mor} in the general
framework of variational analysis,

\begin{Theorem}\label{Theointersecrule} Let $P$ and $\Omega$ be nonempty convex subsets of an LCTV space $X$, where $P$ is a convex polyhedron. Assume that the following quasi-relative interior
qualification  condition \begin{eqnarray}\label{eq1intersecrule} P\cap \qri(\Omega)\ne\emp \end{eqnarray} is satisfied. Then we have the normal cone intersection rule
\begin{equation}\label{eqintersecrule} N(\bar{x};P \cap\Omega)=N(\bar{x};P)+N(\bar{x};\Omega)\; \mbox{\rm for all }\ox\in P\cap \Omega. \end{equation} \end{Theorem} {\bf Proof.}  Fix $\ox\in
P\cap \Omega$ and $\ox^*\in N(\ox; P\cap\Omega)$, and then get by definition \eqref{nc} that \begin{equation*} \la\ox^*,x-\ox\ra\le 0\;\mbox{ for all }\;x\in P\cap\Omega. \end{equation*} Define
further the convex sets in $X\times\R$ by \begin{eqnarray}\label{eqTheta} \begin{array}{ll} &Q:=\big\{(x,\lambda)\in X\times\R\;\big|\;x\in P,\;\lambda\le\la\ox^*,x-\ox\ra\big\}\\ &\mbox{ and }\;
\Theta:=\Omega\times[0,\infty). \end{array} \end{eqnarray} We  have $\qri(\Theta)=\qri(\Omega)\times(0,\infty)$, while $Q$ is a convex polyhedron. It is easy to check from the constructions of the
sets $Q,\Theta$ in \eqref{eqTheta} and the choice of $\ox^*$ that $Q\cap\qri(\Theta)=\emp$. Then the  separation results of Theorem \ref{Sepa-poly-convex} give us a nonzero pair
$(w^*,\gamma)\in X^*\times\mathbb R$ and $\alpha\in \R$ such that \begin{equation}\label{eq2tintersecrule} \la w^*,x\ra+\lambda_1\gamma\le\alpha\leq \la w^*,y\ra+\lambda_2\gamma\;\mbox{ for all }\;(x,\lambda_1)\in
Q,\;(y,\lambda_2)\in\Theta. \end{equation} Furthermore,  there exist pair  $(\tilde{y},\tilde{\lambda}_2)\in\Theta$ satisfying \begin{equation*} \alpha <\la w^*,\tilde{y}\ra+\tilde{\lambda}_2\gamma. \end{equation*}
Using \eqref{eq2tintersecrule} with $(\ox,0)\in Q$ and $(\ox,1)\in\Theta$ yields $\gamma\geq 0$. Now we employ the quasi-relative interior qualification condition \eqref{eq1intersecrule} to show that $\gg\ne 0$. Suppose on the contrary that $\gamma=0$ and then get
\begin{eqnarray*} \begin{array}{ll} &\la w^*,x\ra\le \alpha\leq\la w^*,y\ra\;\mbox{ for all }\;x\in P,\;y\in\Omega,\\ &\mbox{ and }\;\alpha<\la w^*,\tilde{y}\ra\;\mbox{ with }
\tilde{y}\in\Omega. \end{array} \end{eqnarray*} This implies that the sets $P$ and $\Omega$ can be properly separated by a hyperplane $H$ which does not contain $\Omega$, and thus by
Theorem \ref{Sepa-poly-convex} we have $P\cap \qri(\Omega)=\emp$. This contradiction shows that $\gamma>0$.

We immediately deduce from \eqref{eq2tintersecrule} that
\begin{equation*}
\la w^*,\ox\ra\le\alpha\leq\la w^*, x\ra\;\mbox{ for all }\;x\in\Omega,\;\mbox{ and thus }\;-w^*\in N(\ox;\Omega)\;\mbox{ and }\;\disp\frac{-w^*}{\gamma}\in N(\ox;\Omega).
\end{equation*}
It also follows from \eqref{eq2tintersecrule}, when $(x,\beta)\in Q$ with $\beta:=\la \bar{x}^*,x-\ox\ra$ and $(\ox,0)\in\Theta$ that
\begin{equation*}
\la w^*,x\ra+\gamma\la \bar{x}^*,x-\ox\ra\leq \la w^*,\ox\ra\;\mbox{ whenever }\;x\in P.
\end{equation*}
Dividing both sides therein by $\gamma$ gives us the inequality
\begin{equation*}
\Big\la\frac{w^*}{\gamma}+ \bar{x}^*,x-\ox\Big\ra\le 0\;\mbox{ for all }\;x\in P,
\end{equation*}
and so $\frac{w^*}{\gamma}+ \bar{x}^*\in N(\ox; P)$. Therefore we arrive at
\begin{equation*}
\bar{x}^*=(\frac{w^*}{\gamma}+\bar{x}^*)+(\frac{-w^*}{\gamma})\in N(\ox;P)+N(\ox;\Omega),
\end{equation*}
which verifies the inclusion $``\subset$'' in \eqref{eqintersecrule} and thus completes the proof of the theorem since the proof of the opposite inclusion is straightforward.
$\h$

\section{Generalized Differential Calculus under Polyhedrality} \setcounter{equation}{0}

In this section we establish major calculus rules for coderivatives of set-valued mappings and subdifferentials of extended-real-valued functions in convex polyhedral settings of LCTV spaces.

\begin{Definition} \begin{enumerate} \item A set-valued mapping $F\colon X\tto Y$ between LCTV spaces is said to be a {\sc polyhedral set-valued mapping} if $\gph(F)$ is a convex polyhedron in
$X\times Y$. \item A function $f\colon X\to \oR$ defined on an LCTV space $X$ is said to be a {\sc polyhedral function} if $\epi(f)$ is a convex polyhedron in $X\times \R$. \end{enumerate}
\end{Definition}

Recall also that the
(Minkowski) {\em sum} of two set-valued mappings $F_1,F_2\colon X\tto Y$ is \begin{equation*} (F_1+F_2)(x)=F_1(x)+F_2(x):=\big\{y_1+y_2\in Y\;\big|\;y_1\in F_1(x),\;y_2\in F_2(x)\big\},\quad x\in
X. \end{equation*} It is easy to see that $\dom(F_1+F_2)=\dom(F_1)\cap\dom(F_2)$ and that the graph of $F_1+F_2$ is convex provided that both mappings $F_1,F_2$ have this property. Our goal is
to represent the coderivative of the sum $F_1+F_2$ at a given point in terms of the coderivatives of $F_1$ and $F_2$. To proceed, for any $(\ox,\oy)\in\gph(F_1+F_2)$ consider the set
\begin{equation*} S(\ox,\oy):=\big\{(\oy_1,\oy_2)\in Y\times Y\;\big|\;\oy=\oy_1+\oy_2,\;\oy_i\in F_i(\ox)\;\mbox{ as }\;i=1,2\big\} \end{equation*} used in the formulation of the following
coderivative sum rule.

\begin{Theorem}\label{CSR1} Consider two convex set-valued mappings $F_1,F_2\colon X\tto Y$ between LCTV spaces, and let the graphical quasi-relative interior qualification condition
\begin{equation}\label{QCC}
\exists (\Hat x, \Hat{y}_1, \Hat{y}_2)\in X\times Y\times Y\; \mbox{\rm with }(\Hat x, \Hat{y}_1)\in  \gph(F_1) \; \mbox{\rm and }(\Hat x, \Hat{y}_2)\in\qri\big(\gph(F_2)\big)
\end{equation}
be satisfied. Assuming furthermore that  $F_1$ is polyhedral, we have the coderivative sum rule
\begin{equation}\label{csr}
D^*(F_1+F_2)(\ox,\oy)(y^*)=D^*F_1(\ox,\oy_1)(y^*)+D^*F_2(\ox,\oy_2)(y^*)
\end{equation}
valid for all $(\ox,\oy)\in\gph(F_1+F_2)$, for all $y^*\in Y^*$, and for all $(\oy_1,\oy_2)\in S(\ox,\oy)$.
\end{Theorem}
{\bf Proof.} Fix any $x^*\in D^*(F_1+F_2)(\ox,\oy)(y^*)$ and get by equation \eqref{cod} that $(x^*,-y^*)\in N((\ox,\oy);\gph(F_1+F_2))$. For every $(\oy_1,\oy_2)\in S(\ox,\oy)$ consider the convex sets
\begin{eqnarray*}
\begin{array}{ll}
&\Omega_1:=\big\{(x,y_1,y_2)\in X\times Y\times Y\;\big|\;y_1\in F_1(x)\big\},\\
&\Omega_2:=\big\{(x,y_1,y_2)\in X\times Y\times Y\;\big|\;y_2\in F_2(x)\big\}.
\end{array}
\end{eqnarray*}
Then $\Omega_1$ is a convex polyhedron and
\begin{equation*}
\qri(\Omega_2)=\big\{(x,y_1,y_2)\in X\times Y\times Y\;\big|\;(x,y_2)\in\qri(\gph(F_2))\big\}.
\end{equation*}
It follows from the definition that
\begin{equation}\label{cod1}
(x^*,-y^*,-y^*)\in N\big((\ox,\oy_1,\oy_2);\Omega_1\cap\Omega_2\big).
\end{equation}
Furthermore, the qualification condition \eqref{QCC} implies that $\Omega_1\cap\qri(\Omega_2)\ne\emp$. Applying now Theorem~\ref{Theointersecrule} to the set intersection in \eqref{cod1} gives us
\begin{equation*}
(x^*,-y^*,-y^*)\in N\big((\ox,\oy_1,\oy_2);\Omega_1\big)+N\big((\ox,\oy_1,\oy_2);\Omega_2\big).
\end{equation*}
Thus we arrive at the representation
\begin{equation*}
(x^*,-y^*,-y^*)=(x^*_1,-y^*,0)+(x^*_2,0,-y^*)\;\mbox{ with }\;(x^*_i,-y^*)\in N\big((\ox,\oy_i);\gph(F_i)\big)
\end{equation*}
for $i=1,2$. The above representation reads as
\begin{equation*}
x^*=x^*_1+x^*_2\in D^*F_1(\ox,\oy_1)(y^*)+D^*F_2(\ox,\oy_2)(y^*),
\end{equation*}
which justifies the inclusion ``$\subset$" in \eqref{csr}. The opposite inclusion is obvious. $\h$

\begin{Corollary}\label{sr}
Let $f_i\colon X\to\oR$, $i=1,2$, be proper convex functions on an LCTV space $X$, and let
\begin{equation*}\label{riq}
\dom(f_1)\cap\qri\big(\dom(f_2)\big)\ne\emp.
\end{equation*}
Suppose in addition that $f_1$ is polyhedral. Then we have the subdifferential sum rule
\begin{equation}\label{ssr}
\partial(f_1+f_2)(\ox)=\partial f_1(\ox)+\partial f_2(\ox)\; \mbox{\rm for all }\ox\in \dom(f_1)\cap \dom(f_2).
\end{equation}
\end{Corollary}
{\bf Proof.} Define the convex set-valued mappings $F_1,F_2\colon X\tto\R$ by
\begin{equation*}
F_i(x):=\big[f_i(x),\infty\big)\;\mbox{ for }\;i=1,2.
\end{equation*}
Then $\gph(F_1)=\epi(f_1)$ is a convex polyhedron.  Fix any $\Hat x\in\dom(f_1)\cap\qri(\dom(f_2))$ and choose $\gamma>\max\{f_1(\Hat x), f_2(\Hat x)\}$. It is easy to see that $(\Hat x, \gamma)\in \epi(f_1)=\gph(F_1)$ and $(\Hat x, \gamma)\in \qri(\epi(f_2))=\qri(\gph(F_2))$. In addition, for every $x^*\in\partial(f_1+f_2)(\ox)$ we have
\begin{equation*}
x^*\in D^*(F_1+F_2)(\ox,\oy)(1).
\end{equation*}
Applying Theorem~\ref{CSR1} with $\oy_i=f_i(\ox)$ for $i=1,2$ gives us
\begin{equation*}
x^*\in D^*F_1(\ox,\oy_1)(1)+D^*F_2(\ox,\oy_2)(1)=\partial f_1(\ox)+\partial f_2(\ox),
\end{equation*}
which verifies the inclusion ``$\subset$" in \eqref{ssr}. The opposite inclusion is obvious. $\h$

The following lemma was proved in \cite[Theorem 4.1]{CBNG}. We provide here its detailed proof for the convenience of the reader.
\begin{Lemma}\label{svmint} Let $F\colon X\tto Y$ be a convex set-valued mapping between LCTV spaces. Then we have
\begin{equation*}
  \qri(\gph(F)) \supset \{(x,y)\in X\times Y\ \big|\ x\in\qri\big(\dom(F)\big), y\in\sint \big(F(x)\big)\}.
\end{equation*}
\end{Lemma}
{\bf Proof.}  Pick any $(\ox,\oy)\in X\times Y$ with $\ox\in\qri\big(\dom(F)\big)$ and $\oy\in\sint\big(F(\ox)\big)$. To reach a contradiction, suppose that $(\ox,\oy)\notin \qri\big(\gph(F)\big)$. Then
Proposition~\ref{qri-sep1} shows that the sets $\{(\ox,\oy)\}$ and $\gph(F)$ can be properly separated, which means that there exists $(x^*,y^*)\in X^*\times Y^*$ such that
\begin{equation}\label{Eq1-RoctheoGph}
  \la x^*,x\ra +\la y^*,y\ra \leq \la x^*,\ox\ra +\la y^*,\oy\ra \; \text{for all}\; x\in\dom(F), y\in F(x)
\end{equation}
and there exists $(\tilde{x},\tilde{y})\in \gph(F)$ such that
\begin{equation}\label{Eq2-Roctheogph}
  \la x^*,\tilde{x}\ra +\la y^*, \tilde{y}\ra < \la x^*,\ox\ra +\la y^*,\oy\ra.
\end{equation}
Choosing $x=\ox$, \eqref{Eq1-RoctheoGph} shows that
\begin{equation}\label{Eq3-Roctheogph}
  \la y^*,y\ra\leq \la y^*,\oy\ra \; \text{for all} \; y\in F(\ox).
\end{equation} Since $\oy\in\sint\big(F(\ox)\big)$, there exists a symmetric neighborhood $V$ of the origin with $V\subset \sint\big(F(\ox)\big)-\{\oy\}$. It follows from \eqref{Eq3-Roctheogph}
that $\la y^*,v\ra\leq 0$ and $\la y^*,-v\ra \leq 0$ for all $0\ne v\in V$. Thus, $y^*=0$ on $V$. Since $V$ is a symmetric neighborhood of the origin, for any $y\in Y$ there exists $0\ne t\in \R$
such that $ty=v\in V$ and hence $\la y^*,y\ra =\frac{1}{t}\la y^*,v\ra=0$. Thus, $y^*=0$ on $Y$.  It follows from \eqref{Eq1-RoctheoGph} and \eqref{Eq2-Roctheogph} that the sets $\{\ox\}$ and
$\dom(F)$ can be properly separated. Proposition~\ref{qri-sep1} shows that $\ox\notin \qri\big(\dom(F)\big)$. This contradiction completes the proof. $\h$

The theorem below provides sufficient conditions that guarantee the validity of the qualification condition \eqref{QCC}.

\begin{Theorem}\label{CSR2} Let $F_1,F_2\colon X\tto Y$ be two convex set-valued mappings between LCTV spaces. Suppose that the  $F_1$ is polyhedral. Then we have the coderivative sum rule in Theorem \ref{CSR1}  under one of the following conditions:
\begin{enumerate}
  \item $\dom(F_2)$ is quasi-regular and there exists $\Hat u\in\dom(F_1) \cap\qri\big(\dom(F_2)\big)$ such that $\qri\big(F_2(\Hat u)\big)\ne\emptyset$.
  \item There exists  $\Hat u\in\dom(F_1) \cap\qri\big(\dom(F_2)\big)$ such that $\sint\big(F_2(\Hat u)\big)\ne\emptyset$.
  \end{enumerate}
\end{Theorem}
{\bf Proof.} Under assumption (a) there exist $\Hat{z}_1\in F_1(\Hat{u})$ and $\Hat{z}_2\in \qri(F_2(\Hat{u}))$ and thus by \cite[Theorem 4.3(b)]{CBN}  we have
\begin{equation*}
 (\Hat{u},\Hat{z}_1)\in \gph(F_1)\; \text{and}\; (\Hat{u}, \Hat{z}_2)\in\qri\big(\gph(F_2)\big).
\end{equation*}
Thus the qualification condition \eqref{QCC} is satisfied. Similarly, under assumption (b) there exist $\Hat{z}_1\in F_1(\Hat{u})$ and $\Hat{z}_2\in \sint(F_2(\Hat{u}))$. By Lemma  \ref{svmint} the qualification condition \ref{QCC} is also satisfied. Therefore, the conclusion follows directly from Theorem \ref{CSR1}.
$\h$

The next corollary provides another version of the coderivative sum rule replacing the {\em graphical} quasi-relative interior qualification condition \eqref{QCC} with its {\em domain} counterpart. To proceed, we need to mention first some definitions and facts that can be found in \cite{bl}. Recall that a subset $\Omega$ of an LCTV space $X$ is {\em $CS-$closed} if
\begin{equation*}
\Big[\lambda_k\ge 0,\;\sum_{k=1}^\infty\lambda_k=1,\;x_k\in\Omega\;\mbox{ for }\;k\in\N,\;\sum_{k=1}^\infty\lambda_k x_k=x\Big]\Longrightarrow\big[x\in \Omega\big],
\end{equation*}
where the series convergence is with respect to the given topology on $X$. Every $CS-$closed set is clearly convex. If $X$ is Banach, then all the convex sets that are closed, open, finite-dimensional, or $G_\delta$ are $CS-$closed. Recall also that an LCTV space $X$ is {\em Fr\'echet} if it is complete with respect to its topology induced by its translation-invariant metric $d(x,y)$ satisfying $d(x,y)=d(x+z,y+z)$ for all $x,y,z\in X$. Such spaces are the most natural extensions of Banach spaces in the framework of LCTV ones.

\begin{Corollary}{\bf(coderivative sum rule via quasi-relative interiors of domains).}\label{CSR3} Let $F_1,F_2\colon X\tto Y$ be two convex set-valued mappings from an LCTV space $X$ to a separable Fr\'echet space $Y$.
Suppose that $F_1$ is polyhedral, $F_2$ is $CS-$closed and its domain is quasi-regular, and the domain quasi-relative interior qualification condition
\begin{equation}\label{QC1}
\dom(F_1) \cap\qri\big(\dom(F_2)\big)\ne\emp
\end{equation}
is satisfied. Then we have the coderivative sum rule \eqref{csr} valid for all $(\ox,\oy)\in\gph(F_1+F_2)$, for all $y^*\in Y^*$, and for all $(\oy_1,\oy_2)\in S(\ox,\oy)$.
\end{Corollary}
{\bf Proof.} To verify the claimed assertion, it
suffices to show that the domain qualification condition \eqref{QC1} implies the graph one \eqref{QCC}. Having \eqref{QC1}, fix any $\bar{u}\in \dom(F_1)\cap\qri(\dom(F_2))$ and then use the result of \cite[Theorem~2.19]{bl}, which tells us that $F_1(\bar{u})\ne\emptyset$ and $ \qri(F_2(\bar{u}))\ne\emp$  under the assumptions made. Applying Theorem~\ref{CSR2} with condition(a) gives us the result of this corollary. $\h$

Given two set-valued mappings $F\colon X\tto Y$ and $G\colon Y\tto Z$ between LCTV spaces, define their {\em composition} $(G\circ F)\colon X\tto Z$ by
\begin{equation*}
(G\circ F)(x)=\bigcup_{y\in F(x)}G(y):=\big\{z\in G(y)\;\big|\;y\in F(x)\big\},\quad x\in X,
\end{equation*}
and observe that $G\circ F$ is convex provided that both $F$ and $G$ have this property. Fix $\oz\in(G\circ F)(\ox)$ and consider the set
\begin{equation*}
M(\ox,\oz):=F(\ox)\cap G^{-1}(\oz).
\end{equation*}

\begin{Theorem}{\bf(coderivative chain rule via quasi-relative interiors of graphs).}\label{scr} Let $F\colon X\tto Y$ and $G\colon Y\tto Z$ be convex set-valued mappings  between LCTV spaces.
Suppose that one of the two following conditions holds:
\begin{enumerate}
  \item[{\rm (a)}]  $F$ is polyhedral and there exist a triple $(\Hat x,\Hat y,\Hat z)\in X\times Y\times Z$ satisfying
\begin{equation}\label{QC1a}
(\Hat x,\Hat y)\in \gph(F)\;\mbox{\rm and }\;(\Hat y,\Hat z)\in\mbox{\rm qri}\big(\gph(G)\big);
\end{equation}
  \item[{\rm (b)}] $G$ is polyhedral and there exist a triple $(\Hat x,\Hat y,\Hat z)\in X\times Y\times Z$ satisfying
\begin{equation}\label{QC1b}
(\Hat x,\Hat y)\in \mbox{\rm qri}\big(\gph(F)\big)\;\mbox{\rm and }\;(\Hat y,\Hat z)\in \gph(G).
\end{equation}
\end{enumerate}
Then for any $(\ox,\oz)\in\gph(G\circ F)$ and $z^*\in Z^*$ we have the coderivative chain rule
\begin{equation}\label{chain}
D^*(G\circ F)(\ox,\oz)(z^*)=D^*F(\ox,\oy)\circ D^*G(\oy,\oz)(z^*)\;\mbox{ whenever }\;\oy\in M(\ox,\oz).
\end{equation}
\end{Theorem}
{\bf Proof.} Fix $x^*\in D^*(G\circ F)(\ox,\oz)(z^*)$ and $\oy\in M(\ox,\oz)$. Then $(x^*,-z^*)\in N((\ox,\oz);\gph(G\circ F))$ by \eqref{cod}, which means by definition \eqref{nc} that
\begin{equation*}
\la x^*,x-\ox\ra-\la z^*,z-\oz\ra\le 0\;\mbox{ for all }\;(x,z)\in\gph(G\circ F).
\end{equation*}
Define further the two convex subsets of $X\times Y\times Z$ by
\begin{equation*}
\Omega_1:=\gph(F)\times Z\;\mbox{ and }\;\Omega_2:=X\times\gph(G)
\end{equation*}
It easily follows from the constructions above that
\begin{equation*}
(x^*,0,-z^*)\in N\big((\ox,\oy,\oz);\Omega_1\cap\Omega_2\big).
\end{equation*}
The imposed qualification condition (\ref{QC1a}) ensures that $\Omega_1$ is a polyhedral convex set and
\begin{equation*}
  \Omega_1\cap\qri(\Omega_2)\ne\emp,
\end{equation*} while
qualification condition (\ref{QC1b}) ensures that $\Omega_2$ is a polyhedral convex set and
\begin{equation*}
  \qri(\Omega_1)\cap \Omega_2\ne\emp.
\end{equation*}
Then the intersection rule of Theorem~\ref{Theointersecrule} tells us that
\begin{equation*}
(x^*,0,-z^*)\in N\big((\ox,\oy,\oz);\Omega_1\cap\Omega_2\big)=N\big((\ox,\oy,\oz);\Omega_1\big)+N\big((\ox,\oy,\oz);\Omega_2\big).
\end{equation*}
Thus there exists $y^*\in Y^*$ satisfying $(x^*,0,-z^*)=(x^*,-y^*,0)+(0,y^*,-z^*)$ for which we have the normal cone inclusions
\begin{equation*}
(x^*,-y^*)\in N\big((\ox,\oy);\gph(F)\big)\;\mbox{ and }\;(y^*,-z^*)\in N\big((\oy,\oz);\gph(G)\big).
\end{equation*}
This shows by the coderivative definition \eqref{cod} that
\begin{equation*}
x^*\in D^*F(\ox,\oy)(y^*)\;\mbox{ and }\;y^*\in D^*G(\oy,\oz)(z^*),
\end{equation*}
which verifies the inclusion ``$\subset$" in \eqref{chain}. The opposite inclusion is straightforward. $\h$

Similar to the sum rule in Corollary~\ref{CSR2}, we derive from Theorem~\ref{scr} the following coderivative chain rule with the {\em domain} quasi-relative interior qualification condition under mild additional assumptions.

\begin{Corollary}{\bf(coderivative chain rule via quasi-relative interiors of domains).}\label{cod-chain} Let $F\colon X\tto Y$ and $G\colon Y\tto Z$ be convex set-valued mappings between LCTV spaces.
Suppose that one of the two following conditions holds:
\begin{enumerate}
  \item[{\rm (a)}]  $G$ has $CS-$closed values,  $Z$ is a Fr\'echet space,  $F$ is polyhedral, $\dom(G)$ is quasi-regular, and there exists $\Hat x\in \dom(F)$ such that
      \begin{equation*}
        F(\Hat x)\cap \qri\big(\dom(G)\big)\ne\emptyset;
      \end{equation*}
  \item [{\rm (b)}] $G$ is polyhedral, $\dom(F)$ is quasi-regular,  and there exists $\Hat x\in \qri\big(\dom(F)\big)$ such that
      \begin{equation*}
        \qri\big(F(\Hat x)\big)\cap \dom(G)\ne\emptyset.
      \end{equation*}
\end{enumerate}
Then for any $(\ox,\oz)\in\gph(G\circ F)$ and $z^*\in Z^*$ we have the coderivative chain rule \eqref{chain}.
\end{Corollary}
{\bf Proof.}
(a) Fix $\Hat x\in \dom(F)$ and $\Hat y\in F(\Hat x)\cap\qri(\dom(G))$.  Using the imposed additional assumptions tells us by \cite[Theorem~2.19]{bl} that $\qri(G(\Hat y))\ne\emp$. Thus we choose $\Hat z\in\qri(G(\Hat y))$ and derive from \cite[Corollary 4.4]{CBN} that $(\Hat x,\Hat y)\in \gph(F)$ and $(\Hat y,\Hat z)\in\qri(\gph(G))$. This verifies that condition (a) in Theorem~\ref{scr} is satisfied and hence we have \eqref{chain}.

(b) Fix $\Hat x\in \qri\big(\dom(F)\big)$ and $\Hat y\in \qri\big(F(\Hat x)\big)\cap \dom(G)$.  Thus we choose $\Hat z\in G(\Hat y)$ and derive from \cite[Corollary 4.4]{CBN} that $(\Hat x,\Hat y)\in\qri\big( \gph(F)\big)$ and $(\Hat y,\Hat z)\in\gph(G)$. This verifies that condition (b) in Theorem~\ref{scr} is satisfied and hence we have \eqref{chain} here as well.
 $\h$\vspace*{0.05in}

Finally in this section, we present a new subdifferential chain rule for extended-real-valued convex functions on LCTV spaces which follows directly from the coderivative one in Corollary~\ref{cod-chain}. This gives us an infinite-dimensional extension of the classical result of \cite[Theorem~23.9]{r} in finite dimensions without imposing any continuity assumption on the external function that is conventional in the LCTV space framework; see \cite{z}.

\begin{Corollary}{\bf(subdifferential chain rule via quasi-relative interiors of domains).}\label{sub-chain} Let $A\colon X\to Y$ be a linear continuous mapping between LCTV spaces, and let $f\colon Y\to\oR$ be a polyhedral convex function. Assume that the range of $A$ contains a point of $dom(f)$. Then denoting $\oy:=A(\ox)\in\dom(f)$ for some $\ox\in X$ we have the following subdifferential chain rule:
\begin{equation}\label{chain-r1}
\partial(f\circ A)(\ox)=A^*\big(\partial f(\oy)\big):=\big\{A^*y^*\;\big|\;y^*\in\partial f(\oy)\big\}.
\end{equation}
\end{Corollary}
{\bf Proof.} Apply Corollary~\ref{cod-chain} with $F(x):=\{A(x)\}$ and $G(x):=[f(x),\infty)$. Then $\qri(\gph(F))=\gph(A)$, $\dom(G)=\dom(f)$, and $\gph(G)=\epi(f)$. The imposed assumptions guarantee that the qualification condition (b) of Corollary~\ref{cod-chain} is satisfied. This allows us to deduce from
Corollary~\ref{cod-chain} the equalities
\begin{equation*}
\partial(f\circ A)(\ox)=D^*(G\circ A)(\ox,\oy)(1)=D^*A\big(D^*G(\ox,\oy)(1)\big)=A^*\big(\partial f(\oy)\big),
\end{equation*}
which verify \eqref{chain-r1} and thus completes the proof. $\h$

\section{Conjugate Calculus and Duality in Polyhedral Settings}

In this section we first develop refined calculus rules for Fenchel conjugates of convex functions on LCTV spaces with relaxed qualification conditions under appropriate polyhedrality
assumptions. Then we use these results to establish a new duality theorem in a composite format of polyhedral convex optimization. We again employ the geometric approach based on polyhedral
convex separation in LCTV spaces.

Let $\Omega$ be a nonempty subset of a topological vector space $X$. We recall that the support function $\sigma_{\Omega}\colon X^*\to\bar{\R}$ of $\Omega$ is defined by
\begin{equation*}
  \sigma_{\Omega}(x^*):=\sup\{\la x^*,x\ra\;\big|\;x\in\Omega\},\; x^*\in X^*.
\end{equation*}
\begin{Theorem}\label{sigma intersection rule} Let $P$ and  $\Omega$ be nonempty and convex subsets of an LCTV space $X$, where $P$ is a convex polyhedron. Suppose that
\begin{equation*}
    P\cap \qri(\Omega)\neq \emptyset.
\end{equation*}
Then the support function of the intersection $P\cap\Omega$ is represented as \begin{equation}\label{supp1}
\big(\sigma_{P\cap\Omega}\big)(x^*)=\big(\sigma_{P}\s\sigma_{\Omega}\big)(x^*)\;\;\mbox{ for all }\;x^*\in X^*. \end{equation} Furthermore, for any
$x^*\in\dom(\sigma_{P\cap\Omega})$ there exist dual elements $x^*_1,x^*_2\in X^*$ such that $x^*=x^*_1+x^*_2$ and \begin{equation}\label{supp2}
(\sigma_{P\cap\Omega})(x^*)=\sigma_{P}(x^*_1)+\sigma_{\Omega}(x^*_2). \end{equation} \end{Theorem}
{\bf Proof.} The proof of  inequality ``$\le$" in \eqref{supp1} is straightforward.

Now we prove the inequality ``$\ge$" in \eqref{supp1} under the given assumptions. To proceed, fix any $x^*\in\dom(\sigma_{P\cap\Omega})$, denote
$\alpha:=\sigma_{P\cap\Omega}(x^*)\in\R$ for which \begin{equation*} \la x^*,x\ra\le\alpha\;\;\mbox{\rm whenever }\;x\in P\cap\Omega, \end{equation*} and define the nonempty convex subsets of
$X\times\R$ by \begin{eqnarray*}\label{theta1} \begin{array}{ll}  &Q:=\big\{(x,\lambda)\in
X\times\R\;\big|\;x\in P,\;\lambda\le\la x^*,x\ra-\alpha\big\},\\
&\Theta:=\Omega\times[0,\infty).
\end{array} \end{eqnarray*}
It is easy to see from the constructions of $Q$ and $\Th$ that $Q$ is a convex polyhedron with
\begin{equation*}
    Q\cap \qri(\Theta)=\emptyset.
\end{equation*}
Applying Theorem \ref{Sepa-poly-convex} ensures the existence of $(0,0)\ne(w^*,\beta)\in X^*\times\mathbb{R}$ for which \begin{equation}\label{sep f}
\la w^*,x\ra+\beta\lambda_1\le\la w^*,y\ra+\beta\lambda_2\;\mbox{ whenever }\;(x,\lambda_1)\in Q,\;(y,\lambda_2)\in\Theta. \end{equation}
In addition, there are $(\tilde{x},\tilde{\lambda}_1)\in\Theta$ and $(\tilde{y},\tilde{\lambda}_2)\in\Theta$ satisfying
\begin{equation*}
\la w^*,\tilde{x}\ra+\tilde{\lambda}_1\beta<\la w^*,\tilde{y}\ra+\tilde{\lambda}_2\beta.
\end{equation*}
Following the proof of the normal cone intersection rule yields $\bb> 0$.  Taking now the pairs  $(x,\langle x^*,x\rangle-\alpha)\in Q$ in \eqref{sep f} and $(y,0)\in\Theta$, we get  \begin{equation*}  \la
w^*,x\ra+\beta(\la x^*,x\ra-\alpha)\leq \la w^*,y\ra\;\;\mbox{\rm with }\;\bb>0, \end{equation*} which leads us to the estimate \begin{equation*}
\big\la
w^*/\beta+x^*,x\big\ra+\big\la -w^*/\bb,y\big\ra\leq \alpha\;\mbox{ for all }\;x\in P,\;y\in\Omega. \end{equation*} By putting $x^*_1:=w^*/\beta+x^*$ and $x^*_2:=-w^*/\beta$ we arrive at the
inequality ``$\ge$" in \eqref{supp1} and thus get representation \eqref{supp2}. $\h$

The next observation makes a bridge between the conjugate of an arbitrary function and the support function for its epigraph. It is essential in the implementation  of our geometric approach to
conjugate calculus.

\begin{Lemma}\label{Fenchel epi} For any proper function $f\colon X\to\oR$ we have
\begin{equation*}
f^*(x^*)=\sigma_{{\rm\small epi}(f)}(x^*,-1)\;\mbox{ whenever }\;x^*\in X^*.
\end{equation*}
\end{Lemma}
{\bf Proof.}
It follows from the definitions that
\begin{eqnarray*}
\begin{array}{ll}
f^*(x^*)=\sup\big\{\la x^*,x\ra-f(x)\;\big|\;x\in\dom(f)\big\}&=\sup\big\{\la x^*,x\ra-\lambda\;\big|\;(x,\lambda)\in\epi(f)\big\}\\\\
&=\sigma_{{\rm\small epi}(f)}(x^*,-1),
\end{array}
\end{eqnarray*}
which therefore verifies the claimed formula. $\h$

Here is the {\em conjugate sum rule} in LCTV spaces.

\begin{Theorem}\label{Fenchel sum rule} Let $f,g\colon X\to\oR$ be proper convex functions on an LCTV space $X$, where $f$ is a polyhedral function. Suppose that
\begin{equation*}
    \dom(f)\cap \qri(\dom(g))\neq \emptyset.
\end{equation*}
Then we have the conjugate sum rule \begin{equation}\label{Fenchelsum} (f+g)^*(x^*)=\big(f^*\s g^*\big)(x^*)\;\mbox{ for all }\;x^*\in X^*. \end{equation} Moreover, the infimum in $(f^*\s g^*)(x^*)$
is attained, i.e., for any $x^*\in\dom(f+g)^*$ there exist vectors $x^*_1,x^*_2\in X^*$ for which \begin{eqnarray*}\label{inf-conj} (f+g)^*(x^*)=f^*(x^*_1)+g^*(x^*_2),\quad x^*=x^*_1+x^*_2.
\end{eqnarray*} \end{Theorem} {\bf Proof.} Fixing any $x^*_1,x^*_2\in X^*$ with $x^*_1+x^*_2=x^*$, we get \begin{eqnarray*} \begin{array}{ll} f^*(x^*_1)+g^*(x^*_2)&=\sup\big\{\la
x^*_1,x\ra-f(x)\;\big|\;x\in X\big\}+\sup\big\{\la x^*_2,x\ra-g(x)\;\big|\;x\in X\big\}\\ &\ge\sup\big\{\la x^*_1,x\ra-f(x)+\la x^*_2,x\ra-g(x)\;\big|\;x\in X\big\}\\ &=\sup\big\{\la
x^*,x\ra-(f+g)(x)\;\big|\;x\in X\big\}=(f+g)^*(x^*). \end{array} \end{eqnarray*} Note that this inequality does not require the continuity assumption.

Let us prove that $(f^*\s g^*)(x^*)\le(f+g)^*(x^*)$. We only need to consider the case where $(f+g)^*(x^*)<\infty$. Define two convex sets by \begin{eqnarray}\label{omega-conj}
\begin{array}{ll} &P:=\big\{(x,\lambda_1,\lambda_2)\in X\times\R\times\R\;\big|\;\lambda_1\ge f(x)\big\}=\epi(f)\times\R,\\ &\Omega:=\big\{(x,\lambda_1,\lambda_2)\in
X\times\R\times\R\;\big|\;\lambda_2\ge g(x)\big\}. \end{array} \end{eqnarray} Similarly to Lemma~\ref{Fenchel epi} we get the representation \begin{equation}\label{conj-supp}
(f+g)^*(x^*)=\sigma_{P\cap\Omega}(x^*,-1,-1). \end{equation}
We can show that $P$ and $\Omega$ satisfy the support function intersection rule requirements. Applying this theorem to the right-hand side of \eqref{conj-supp} gives us triples $(x_1^*,-\alpha_1,-\alpha_2)\in
X^*\times\R\times\R$ and $(x_2^*,-\beta_1,-\beta_2)\in X^*\times\R\times\R$ such that $(x^*,-1,-1)=(x_1^*,-\alpha_1,-\alpha_2)+(x_2^*,-\beta_1,-\beta_2)$ and \begin{equation*}
(f+g)^*(x^*)=\sigma_{P\cap\Omega}(x^*,-1,-1)=\sigma_{P}(x^*_1,-\alpha_1,-\alpha_2)+\sigma_{\Omega}(x^*_2,-\beta_1,-\beta_2). \end{equation*} If $\alpha_2\ne 0$, then
$\sigma_{P}(x_1^*,-\alpha_1,-\alpha_2)=\infty$, which is not possible. Thus $\alpha_2=0$ and similarly $\beta_1=0$. Employing now Lemma~\ref{Fenchel epi} and taking into account the structures
of the sets $P$ and $\Omega$ in \eqref{omega-conj} imply that \begin{eqnarray*} \begin{array}{ll}
(f+g)^*(x^*)&=\sigma_{P\cap\Omega}(x^*,-1,-1)=\sigma_{P}(x^*_1,-1,0)+\sigma_{\Omega}(x^*_2,0,-1)\\ &=\sigma_{{\rm\small epi}(f)}(x^*,-1)+\sigma_{{\rm\small epi}(g)}(x^*_2,-1)\\
&=f^*(x^*_1)+g^*(x^*_2)\ge\big(f^*\s g^*\big)(x^*). \end{array} \end{eqnarray*} This justifies the sum rule \eqref{Fenchelsum} and that the infimum is attained. $\h$

\begin{Theorem}\label{Theo_Fen_Poly} Let $f,g \colon X\to \oR$ be convex functions. Suppose that $f$ is polyhedral and that
\begin{equation*}
  \dom(f)\cap \qri\big(\dom(g)\big)\ne\emptyset.
\end{equation*}
Then
\begin{equation*}\label{duality_Theo1}
    \inf\{f(x)+g(x)\ |\ x\in X\}=\sup\{-f^*(-x^*)-g^*(x^*)\ |\ x^*\in X^*\}
\end{equation*}
\end{Theorem}
{\bf Proof.}
The inequality $``\geq$" follows directly from the definition.

We next prove the inequality $``\leq$". Since the this inequality is obvious when
\begin{equation*}
    \inf\{f(x)+g(x)\; \big|\; x\in X\}=-\infty,
\end{equation*}
it suffices to consider the case where $\alpha:=\inf\{f(x)+g(x)\; \big|\; x\in X\}\in \R$. It is easy to see that
\begin{equation*}
    \alpha=\inf\{f(x)+g(x)\; \big|\; x\in X\}=-\sup\{\la 0,x\ra-(f+g)(x)\; \big|\; x\in X\}=-(f+g)^*(0).
\end{equation*}
Using conjugate sum rule from Theorem~\ref{Fenchel sum rule}, we can find $x^*\in X^*$ such that
\begin{equation*}
    \alpha=-(f+g)^*(0)=-f^*(-x^*)-g^*(x^*).
\end{equation*}
Therefore, we arrive the inequality $``\leq$" and hence the Theorem has been proved.
 $\h$

We now investigate the {\em composite convex optimization} framework of Fenchel duality in general LCTV spaces with specifying the results in Banach and finite-dimensional
settings. Given proper convex functions $f\colon X\to\oR$,  $g\colon Y\to\oR$ and a linear continuous operator $A\colon X\to Y$ between LCTV spaces $X$ and $Y$, consider the following {\rm primal}
minimization problem: \begin{eqnarray}\label{FP}\mbox{ minimize }&f(x)+g(Ax)\;\mbox{ subject to }\;x\in X.\end{eqnarray} Note that the optimization problem \eqref{FP} is written in the
unconstrained format, while it actually incorporates constraints via the domains of the extended-real-valued functions $f$ and $g$. Using the conjugate functions of $f,g$ and the adjoint operator
of $A$, the {\rm Fenchel dual problem} of \eqref{FP} is defined in the maximization form as follows: \begin{eqnarray}\label{FD} \mbox{ maximize }&-f^*(A^*y^*)-g^*(-y^*)\;\mbox{ subject to
}\;y^*\in Y^*. \end{eqnarray} As we know, the conjugate functions $f^*\colon X^*\to\oR$ and $g^*\colon Y^*\to\oR$ are convex while entering \eqref{FD} with the negative sign, and thus the Fenchel
dual problem also belongs to convex minimization.\vspace*{0.03in}

The following result establishes an inequality relationship between optimal values of the problem and dual problems. This relationship, which is sometimes called {\em weak duality}, holds with no
convexity assumptions of $f$ and/or $g$, and its proof is a consequence of the definitions.
\begin{Theorem}\label{Fenchel chain rule}
Let $A\colon X\to Y$ be a continuous linear mapping between LCTV spaces, and let $g\colon Y\to\overline{\R}$ be a proper convex function. Suppose that  $g$ is polyhedral and
\begin{equation}\label{eq1_composi_fen}
    AX \cap \dom g \ne\emptyset.
\end{equation}
Then we have the conjugate chain rule
\begin{equation*}
    (g\circ A)^*(x^*)=\inf\{g^*(y^*)\; \big|\; y^*\in (A^*)^{-1}(x^*)\},\; x^*\in X^*.
\end{equation*}
Furthermore, for any $x^*\in\dom(g\circ A)^*$ there exists $y^*\in (A^*)^{-1}(x^*)$ such that
\begin{equation*}
    (g\circ A)^*(x^*)=g^*(y^*).
\end{equation*}
\end{Theorem}
{\bf Proof.} Picking $y^*\in(A^*)^{-1}(x^*)$ gives us by definition that \begin{eqnarray*} \begin{array}{ll}
g^*(y^*)&=\sup\big\{\la y^*,y\ra-g(y)\;\big|\;y\in Y\big\}\\ &\ge\sup\big\{\la y^*,Ax\ra-g(Ax)\;\big|\;x\in X\big\}\\ &=\sup\big\{\la A^*y^*,x\ra-(g\circ A)(x)\;\big|\;x\in X\big\}\\
&=\sup\big\{\la x^*,x\ra-(g\circ A)(x)\;\big|\;x\in X\big\}=(g\circ A)^*(x^*). \end{array} \end{eqnarray*} This implies the inequality \begin{equation*}
\inf\big\{g^*(y^*)\;\big|\;y^*\in(A^*)^{-1}(x^*)\big\}\ge(g\circ A)^*(x^*). \end{equation*}
Note that this inequalities holds obviously if $(A^*)^{-1}(x^*)=\emptyset$.

Let us show that the opposite inequality is also satisfied. We can assume that $x^*\in\dom((g\circ
A)^*)$ and then define the convex sets
\begin{equation}\label{omega-chain}\;P:=X\times\epi(g)\subset X\times Y\times\R\mbox{ and }\;\Omega:=\gph(A)\times\R. \end{equation}
It follows directly from the above constructions that
\begin{equation*} (g\circ A)^*(x^*)=\sigma_{P\cap \Omega}(x^*,0,-1)<\infty. \end{equation*}
We can show that $\qri(\Omega)=\Omega$. It follows from the assumption $g$ is polyhedral that $P$ is also polyhedral. Using identity \eqref{eq1_composi_fen} gives us the qualification
\begin{equation*}
    P \cap \qri(\Omega)\ne\emptyset.
\end{equation*}
Then Theorem~\ref{sigma intersection rule} tells us that there exist triples $(x^*_1,y^*_1,\alpha_1)$ and
$(x^*_2,y^*_2,\alpha_2)$ in the space $X^*\times Y^*\times\R$ satisfying \begin{equation*} (x^*,0,-1)=(x^*_1,y^*_1,\alpha_1)+(x^*_2,y^*_2,\alpha_2)\;\mbox{ and} \end{equation*} \begin{equation*}
\sigma_{P\cap\Omega}(x^*,0,-1)=\sigma_{P}(x^*_1,y^*_1,\alpha_1)+\sigma_{\Omega}(x^*_2,y^*_2,\alpha_2). \end{equation*} It follows from the structures of $P$ and  $\Omega$  in
\eqref{omega-chain} that $\alpha_1=0$ and $x^*_2=0$. This gives us the representation \begin{equation*}
\sigma_{P\cap\Omega}(x^*,0,-1)=\sigma_{P}(x^*,-y^*_2,0)+\sigma_{\Omega}(0,y^*_2,-1) \end{equation*} for some $y^*_2\in Y^*$. Thus we arrive at the equalities \begin{eqnarray*}
\begin{array}{ll} \sigma_{P\cap\Omega}(x^*,0,-1)&=\sup\big\{\la x^*,x\ra-\la y^*_2,Ax\ra\;\big|\;x\in X\ra\big\}+\sigma_{{\rm\small epi}(g)}(y^*_2,-1)\\ &=\sup\big\{\la
x^*-A^*y_2^*,x\ra\;\big|\;x\in X\big\}+g^*(y^*_2), \end{array} \end{eqnarray*} which allow us to conclude that $x^*=A^*y_2^*$ and therefore \begin{equation*}
\sigma_{P\cap\Omega}(x^*,0,-1)=g^*(y^*_2)\ge\inf\big\{g^*(y^*)\;\big|\;y^*\in(A^*)^{-1}(x^*)\big\}. \end{equation*} This justifies both statements of the theorem.
$\h$

\begin{Proposition}\label{wd} Consider the optimization problem \eqref{FP} and its dual \eqref{FD} in LCTV spaces, where the functions $f$ and $g$ are not assumed to be convex. Define the optimal
values in these problems by \begin{eqnarray*} \begin{array}{ll} &\disp\Hat p:=\inf_{x\in X}\big\{f(x)+g(Ax)\big\},\\ &\disp\Hat d:=\sup_{y^*\in Y^*}\big\{-f^*(A^*y^*)-g^*(-y^*)\big\}. \end{array}
\end{eqnarray*} Then we always have $\Hat p\ge\Hat d$. \end{Proposition} {\bf Proof.} Directly from the definitions of conjugate functions and adjoint operators, for any $y^*\in Y^*$ we have
\begin{align*} -f^*(A^*y^*)-g^*(-y^*)&=-\sup_{x\in X}\{\la A^*y^*,x\ra-f(x)\}-\sup_{y\in Y}\{\la-y^*,y\ra-g(y)\}\\ &=\inf_{x\in X}\{f(x)-\la y^*,Ax\ra\}+\inf_{y\in Y}\{g(y)+\la y^*,y\ra\}\\
&\leq\inf_{x\in X}\{f(x)-\la y^*,Ax\ra\}+\inf_{x\in X}\{g(Ax)+\la y^*,Ax\ra\}\\ &\leq\inf_{x\in X}\{f(x)-\la y^*,Ax\ra+g(Ax)+\la y^*,Ax\ra\}\\ &=\inf_{x\in X}\{f(x)+g(Ax)\}=\Hat p. \end{align*}
Taking the supremum with respect to all $y^*\in Y^*$ yields $\Hat d\leq\Hat p$. $\h$

\begin{Theorem}\label{sd} Consider the optimization problem \eqref{FP} and its dual \eqref{FD}. In addition to the standing assumptions on $f,g$, and $A$ in the formulation of the problem,
suppose that $g$ is polyhedral and
    \begin{equation}\label{QCD}
        \dom(g\circ A)\cap \qri\big(\dom (f)\big)\ne\emptyset.
    \end{equation}
 Then we have the equality $\Hat p=\Hat d$. Furthermore, if the number $\Hat p$ is finite, then the
supremum in the definition of $\Hat d$ is attained. \end{Theorem}
{\bf Proof.}
Due to Proposition~\ref{wd}, it remains to show that $\Hat p\le\Hat d$. Since the latter inequality is obvious when
$\Hat p=-\infty$, it suffices to consider the case where $\Hat p\in\R$. We clearly have the equalities \begin{eqnarray*}\begin{array}{ll}\Hat p&:=\disp\inf_{x\in X}\big\{f(x)+
g(Ax)\big\}=-\disp\sup_{x\in X}\big\{\la 0,x\ra-[f+(g\circ A)](x)\big\}\\&=-[f+(g\circ A)]^*(0). \end{array}\end{eqnarray*}
Since $g$ is a convex polyhedral function, using \cite[Theorem 3.2]{luan2017} we can show $g\circ A$ is also a convex polyhedral function. By the conjugate sum rule from Theorem~\ref{Fenchel sum rule}, we find $x^*\in X^*$ such that \begin{equation*} \Hat p=-[f+(g\circ A)]^*(0)=-f^*(-x^*)-(g\circ A)^*(x^*). \end{equation*} The
conjugate chain rule from Theorem~\ref{Fenchel chain rule} gives us $y^*\in Y^*$ satisfying \begin{equation*} A^*y^*=x^*,\;\;\mbox{ and }\;(g\circ A)^*(x^*)=g^*(y^*). \end{equation*} Therefore,
we arrive at the relationships \begin{equation*}\Hat p=-f^*(-A^*y^*)-g^*(y^*)\le\Hat d,\end{equation*} which also ensure that the supremum in the definition of $\Hat d$ is attained.
$\h$

\begin{Corollary} The qualification condition \eqref{QCD} is satisfied if
\begin{equation*}
    \dom(g)\cap A(\qri(\dom(f))\neq\emptyset.
\end{equation*}
\end{Corollary}
{\bf Proof.} Choose $y_0\in \dom(g)$ such that $y_0=A(x_0)$ for some $x_0\in \qri(\dom(f))$. Then $x_0\in \dom(g\circ A)\cap \qri(\dom(f))$. $\h$

 \end{document}